\documentclass{amsproc}
\usepackage{amsfonts,eucal}
\input xy
\xyoption{all}
\begin{document}
\newcommand{\A}{{\mathcal A}}
\newcommand{\BU}{{B{\rm U}}}
\newcommand{\C}{{\mathbb C}}
\newcommand{\Con}{{\mathcal C}}
\newcommand{\Dim}{{\rm dim}}
\newcommand{\m}{{\mathcal M}}
\newcommand{\M}{{\overline{\mathcal M}}}
\newcommand{\Q}{{\mathbb Q}}
\newcommand{\R}{{\mathbb R}}
\newcommand{\Z}{{\mathbb Z}}
\newcommand{\Diff}{{\rm Diff}}
\newcommand{\Metrics}{{\rm Metrics}}
\newcommand{\MU}{{\bf MU}}
\newcommand{\la}{{\langle}}
\newcommand{\ra}{{\rangle}}
\newcommand{\proto}{{\mathcal Smpl}}
\newcommand{\AC}{{\mathcal Acx}}
\newcommand{\el}{{\mathcal L}}
\newcommand{\W}{{\mathbb W}}
\newcommand{\half}{{\textstyle {\frac{1}{2}}}}

\title {Topological gravity in dimensions two and four}\bigskip
\author{Jack Morava}
\address{Department of Mathematics, Johns Hopkins University,
Baltimore, Maryland 21218}
\email{jack@math.jhu.edu}
\thanks{The author was supported in part by the NSF}
\date {21 June 1999}
\begin{abstract} Recent work on gravity in two dimensions generalizes
naturally to four dimensions. This is a version of a talk at the 1999 
operads conference in Utrecht; the bulk of the paper [\S 2] is 
examples. \end{abstract}

\maketitle

\section{Basic definitions}

\noindent
{\bf 1.1} The (symmetric monoidal) {\bf two}-category 
$$(Gravity)_{d+1}$$ has {\bf objects:} 
compact oriented $d$-manifolds, with \medskip

\noindent
$\bullet$ {\bf morphisms} $V_0 \rightarrow V_1$ : $(d+1)$-manifolds 
$W$ with $\partial W \cong V_0^{op} \sqcup V_1$, and \medskip

\noindent
$\bullet$ diffeomorphisms $\tilde W \rightarrow W$ as {\bf two-morphisms}.
\medskip

\noindent
The category $Mor(V_0,V_1)$ with cobordisms from $V_0$ to $V_1$ as
objects and diffeomorphisms (equal to the identity on the
boundary) as morphisms, is a hom-object in this two-category. Disjoint
union defines the {\bf monoidal structure}, and the category has an 
orientation-reversing {\bf adjoint equivalence} with its opposite.
\medskip

\noindent
{\bf 1.2} The {\bf topological} category $$({\rm Gravity})_{d+1}$$ 
has compact {\bf Riemannian} $d$-manifolds as objects, and the spaces
$${\rm Mor}(V_0,V_1) := \coprod_{V_0^{op} \sqcup V_1 \cong \partial W} 
(\Metrics/\Diff)(W)$$ as its hom-objects. Alternately: a morphism is a
$(d+1)$-dimensional cobordism, together with (the equivalence
class of) a Riemannian metric on it. \medskip

\noindent
The group of diffeomorphisms which fix a frame at a point acts
{\bf freely} on the space of Riemannian metrics on a complete
manifold, so the morphism spaces of $({\rm Gravity})_{d+1}$
are roughly just the classifying spaces [12] of the morphism categories of
$(Gravity)_{d+1}$. \medskip 

\noindent
{\bf 1.3} When $d=1$ the group of diffeomorphisms of a Riemann surface (of
genus $>1$) has {\bf contractible} components, and its mapping-class
group $\Gamma = \pi_0(\Diff)$ acts with {\bf finite} isotropy
on Teichm\"uller space, defining a {\bf rational} homology isomorphism 
$$B\Diff \sim B\Gamma \sim {\rm Teich} \times_{\Gamma} E\Gamma 
\rightarrow {\rm Teich} \times_{\Gamma} {\rm pt} = {\mathcal M}\;.$$ 
$(\rm Gravity)_{1+1}$ is thus very similar to the category Segal
constructed to define conformal field theory. \medskip 

\noindent
{\bf 1.4} A monoidal functor from the topological gravity
category to some simpler monoidal category, such as Hilbert spaces 
and trace-class maps, or modules over a ring spectrum, defines a 
{\bf theory of topological gravity}. Replacing the Hom-objects in
this category with their sets of components leads to TFT's in the 
sense of Atiyah, and replacing the Hom-objects with their rational 
homotopy types leads to cohomological field theories.
\bigskip
 
\section {Some examples}

\noindent
{\bf 2.1} The `intersection homology' of a connected surface 
$$\Sigma \mapsto \ker [H^{1}(\Sigma,\C) \rightarrow H^{1}(\partial 
\Sigma,\C)]$$ is a simple example: if $\Sigma \circ \Sigma'$ is 
the composition of two surfaces along a boundary component, 
then the induced map $$\tau : B\Diff(\Sigma) \rightarrow \BU$$ 
fits in the commutative diagram 
$$\xymatrix{
{B\Diff(\Sigma) \times B\Diff(\Sigma')} \ar[r] \ar[d]^{\tau \times \tau} & 
{B\Diff(\Sigma \circ \Sigma')} \ar[d]^{\tau}  \\
{\BU \times \BU} \ar[r]^{\oplus} & {\BU} \;.}
$$ 
This defines a theory of topological gravity with values in a monoidal 
topological category having {\bf one} object, with the $H$-space $\BU$ of 
morphisms. [If we want to fuse along more than one boundary component,
though, we need to be more careful.]
\medskip

\noindent
This functor has more structure: it takes values in symplectic
lattices. The composition $$B{\rm Sp}(\Z) \rightarrow B{\rm Sp}(\R) \sim
\BU \rightarrow B({\rm U/Sp}) \sim {\rm SO/U}$$ is a rational homology 
isomorphism, so $\tau$ lifts to a map $$\tau_{T} : B\Diff(\Sigma)
\rightarrow {\rm SO/U}$$ which sends a surface to its harmonic one-forms
with the complex structure defined by the Hodge $*$-operator. This is
a form of the Abel-Jacobi-Torelli map $$\Sigma \mapsto 
{\rm SP}^{\infty} \Sigma : {\mathcal M} \rightarrow {\mathcal A}$$ 
which sends a surface to its Jacobian; note that the union of $\Sigma$
with its opposite defines a quaternionic object, which $\tau_T$
maps to zero. \medskip

\noindent
{\bf 2.2} Kontsevich-Witten theory is a much deeper example, with the 
(rationalized) complex cobordism ring-spectrum as target object. 
A toy version is easy to construct: \medskip

\noindent
Suppose $\partial \Sigma$ has at most one component; then capping it off
defines the closed surface $\Sigma_{D} := \Sigma \circ D$. 
The cobordism class of the bundle $$[\Sigma_{D}]: \Sigma_{D} 
\times_{\Diff(\Sigma)} E\Diff(\Sigma) \rightarrow B\Diff(\Sigma)$$ is 
{\bf primitive}, in the sense that it behaves additively under composition: 
the pullback $$\mu^{*}[\Sigma \circ \Sigma'] \in MU^{-2}(B\Diff(\Sigma) 
\times B\Diff(\Sigma'))$$ under the composition
$$\mu: \Diff(\Sigma) \times \Diff(\Sigma') \rightarrow \Diff(\Sigma 
\circ \Sigma')$$ is the sum $$\epsilon^{*} [\Sigma_{D}] 
\otimes 1 + 1 \otimes {\epsilon'}^{*} [\Sigma'_{D}]\;,$$
where $$\epsilon : \Diff(\Sigma) \rightarrow \Diff(\Sigma_{D})$$ 
extends the diffeomorphism by the identity [9]. \medskip

\noindent
This says we can pull $\Sigma \circ \Sigma'$ apart as if it were 
made of taffy: the standard family of quadratic cones in $\R^3$ glued to 
$(\Sigma^{op} \sqcup \Sigma') \times I$ defines a $\Diff(\Sigma) 
\times \Diff(\Sigma')$-equivariant cobordism from $\Sigma \circ \Sigma'$ 
to $\Sigma_{D} \sqcup \Sigma'_{D}$. \medskip

\noindent
If we tensor with $\Q$ [and supress the grading] then the map $\tau_{kw}$
representing $$\exp([\Sigma]) \in MU_{\Q}^{*}(B\Diff(\Sigma))$$ fits in the 
homotopy-commutative diagram 
$$\xymatrix{
{B\Diff(\Sigma)^+ \wedge B\Diff(\Sigma')^+} \ar[r] \ar[d]^{\tau_{kw} 
\wedge \tau_{kw}} & {B\Diff(\Sigma \circ \Sigma')^+} \ar[d]^{\tau_{kw}}  \\
{\MU_{\Q} \wedge_{\Q} \MU_{\Q}} \ar[r] & {\MU_{\Q}\;.}}
$$ In both these examples, a surface is sent to some kind of 
configuration space of points: such constructions take unions 
to products. Proper Kontsevich-Witten theory involves
much more complicated configuration spaces. \medskip

\noindent
These constructions capture much of what's known about the stable
cohomology of moduli spaces. Mumford's conjecture, for example, is
equivalent to the assertion that $\tau_{kw}$ defines an isomorphism
on rational cohomology. \medskip

\noindent
{\bf 2.3} The {\bf Floer homology} $HF^*(Y)$ of a compact 3-manifold $Y$
(e.g. a homology sphere) is defined by the Chern-Simons
functional on the space of connections $\Con (Y)$ mod gauge
equivalence on a (trivial) $G$-bundle over $Y$. It is
periodically graded, and has a kind of Poincar\'e pairing.
Following Cohen, Jones, and Segal [3], I assume it is defined by an
underlying spectrum ${\bf HF}(Y)$. Because the
space of connections satisfies $$\Con (Y_0 \sqcup Y_1) = \Con (Y_0)
\times \Con (Y_1) \;,$$ it follows that $${\bf HF}(Y_0 \sqcup Y_1)
= {\bf HF}(Y_0) \wedge {\bf HF}(Y_1) \;.$$ Atiyah [1] saw that Floer 
homology is a topological field theory: when $Y$ bounds $Z$, the space 
$\A(Z)$ of Yang-Mills instantons on $Z$ defines (by restriction to 
$\partial Z$) a kind of Lagrangian cycle  $$[\A(Z) \rightarrow \Con (Y)] 
\in HF^*(Y) \;,$$ and if $$\partial Z = Y^{op} \sqcup Y', \; \partial Z'
 = {Y'}^{op} \sqcup Y''$$ then $[\A(Z)] \wedge [\A(Z')]$ should map 
to $[\A(Z \cup_{Y} Z']$ under the pairing $${\bf HF}(Y^{op} \sqcup Y') 
\wedge {\bf HF}({Y'}^{op} \sqcup Y'') \rightarrow {\bf HF}(Y^{op} 
\sqcup Y'') \;.$$ In fact Yang-Mills on $Z$ presupposes a Riemannian metric, 
and there is a {\bf family} $$\A(Z) \times_{\Diff(Z)} E\Diff(Z) 
\rightarrow \Con (\partial Z) \times B\Diff(Z)$$ of Lagrangian cycles;
its hypothetical class $$\tau_A : B\Diff(Z)^+ \rightarrow 
{\bf HF}(\partial Z)$$ should define a theory of topological gravity.
\medskip

\noindent
{\bf 2.4} In honest Kontsevich-Witten theory [6,13] the analogue 
$\tau_{\bf kw}$ of $\exp([\Sigma])$ is the class 
$$\sum_{n \geq 0} \la \M_g^n \ra /n! \in MU_{\Q}^*(\M_g)\;,$$ 
where $\la \M_g^n \ra$ is the cobordism 
class of a forgetful map to the Deligne-Mumford space of stable 
algebraic curves of genus $g$, from
a compactification of the space of smooth curves marked with $n$ distinct
points. Its characteristic number polynomial is $$\Phi_*{\bf m}_{tot}
(-\nu_{fake}) \in H^*(\M_g,MU_{\Q}) \;,$$ where $\Phi$ is the forgetful
map from the Deligne-Mumford-Knudsen space $\M_g^n, \; {\bf m}_{tot}$ is
the characteristic class defined by the total monomial 
symmetric function, and $-\nu_{fake}$ is the sum of the tangent line 
bundles to the modular curve at its marked points. I'm indebted to Gorbounov, 
Manin, and Zograf for correcting my mistaken assertion [10] that $\nu_{fake}$ 
is the formal normal bundle of $\Phi$ : above
the divisor $\M_g^{n-k+1} \times \M_0^{1+k}$ on $\M_g^n$ defined by curves 
with two irreducible components, one of genus zero, these two bundles differ
stably by the sum of the pair of tangent lines at the double point. \medskip

\noindent
I claim that $\tau_{\bf kw}$ respects a monoidal structure defined
by Knudsen's gluing map $\mu$: in the simplest case this means that
$$\xymatrix{
{\M_g^{1+} \wedge \M_h^{1+}} \ar[r] \ar[d] & 
{\M_{g+h}^+} \ar[d]  \\
{\MU_{\Q} \wedge_{\Q} \MU_{\Q}} \ar[r] & {\MU_{\Q}}}$$ 
commutes, or equivalently that $$\mu^* \la \M_{g+h}^n \ra = \sum_{p+q=n} \la 
\M_g^{p:1} \ra \times \la \M_h^{q:1} \ra \;,$$ where $\la \M_g^{p:r} \ra \in
MU_{\Q}^*(\M_g^r)$ is defined by the partially forgetful morphism $\Phi_{r} 
\;:\; \M_g^{p+r} \rightarrow \M_g^p \;.$ This follows because the diagram 
$$\xymatrix{
{\coprod_{p+q=n} \M_g^{p+1} \times \M_h^{1+q}} \ar[r] \ar[d]^{\tilde \Phi_1
\times \tilde \Phi_1} & {\M_{g+h}^n} \ar[d]^{\Phi}  \\ {\M_g^1 \times \M_h^1}
\ar[r]^{\mu} & {\M_{g+h}} }$$ 
is a pullback in the category of smooth stacks, so 
$$\mu^* \Phi_* = (\tilde \Phi_1 \times \tilde \Phi_1)_* 
\tilde \mu^*$$ by the base-change theorem [8] of Moerdijk and Pronk.
Now $\tilde \mu^*$ pulls back tangent lines at 
marked points to tangent lines at marked points, so $$ \mu^*_{p,q} 
(\nu_{fake}^{(n)}) = \nu_{fake}^{(p)} \otimes 1 + 1 \otimes 
\nu_{fake}^{(q)}\;,$$ and ${\bf m}_{tot}$ is multiplicative, so 
$$\mu^* \la \M_{g+h}^n \ra = 
\mu^* \Phi_* {\bf m}_{tot}(-\nu_{fake}^{(n)})$$ is a sum of terms of the 
form $$(\tilde \Phi_1\times \tilde \Phi_1)_*({\bf m}_{tot}(-\nu_{fake}^{(p)})
\times {\bf m}_{tot}(-\nu_{fake}^{(q)})) \;,$$ {\bf QED.} \medskip

\noindent
{\bf 2.5i} Topological gravity coupled to the quantum cohomology of a 
smooth projective variety $V$ is a (largely) conjectural theory defined, 
in the terms suggested here,  by $$\sum_{n \geq 0} 
\la \M_g^{n:k}(V) \ra /n! \in MU_{\Q}^*(\M_g^{k+} \wedge V^{k+}) \;.$$ 
The configuration spaces are now suitable moduli spaces of stable 
maps [2] from curves to $V$, graded by degree (in $H_2(V,\Z)$); but this 
grading will be supressed here. The representing morphism 
$$\tau_{\bf kw}(V) : \M_g^{k+} \rightarrow [V^{k+},\MU_{\Q}]$$ hypothetically 
defines a functor to a monoidal category with objects $\mathbb N$, morphisms
$${\rm mor}(j,k) = [V^{j+k+},\MU_{\Q}] \;,$$ and compositions defined
by a Poincar\'e trace $$ [V^+ \wedge V^+,\MU_{\Q}] \rightarrow 
\MU_{\Q} \;;$$ the superscript $+$'s indicate the addition of a disjoint
basepoint, according to the conventions of homotopy theory. There is a 
natural monoidal functor $n \mapsto [V^{n+},\MU_\Q]$ 
to the usual category of $\MU_\Q$-module spectra; this is the first 
two-dimensional case in which the target category for a topological gravity
theory has not been slightly degenerate. \medskip

\noindent
The simplest case of the monoidal axiom asserts the commutativity of
$$\xymatrix{
{\M_g^{1+} \wedge \M_h^{1+}} \ar[r]^{\mu} \ar[d]  & {\M_{g+h}^+} \ar[d]  \\
{ [V,\MU_{\Q}] \wedge_{\MU_{\Q}} [V,\MU_{\Q}}] \ar[r] & 
{\MU_{\Q}}\;;}$$
the map from $\M_0^{3+}$ to $ [V^{3+},\MU_{\Q}]$ then defines a
new (quantum) product. When $V$ is a point, this agrees with the usual 
product on $\MU_{\Q}$, after a renormalization which shifts the identity
element: the simplified model of \S 2.2 multiplies the identity by $q = 
\exp (\C P_1))$, while in proper Kontsevich-Witten theory the analogous
$q$ is a generating function for the cobordism classes of the spaces 
$\M_0^{n+3}$. Here the word renormalization
is meant quite literally: $q$ accounts for the cloud of virtual particles
familiar in other quantum contexts. \medskip

\noindent
A relative version of this construction uses 
Gromov-Witten classes $$\la \widetilde{\mathcal M}_g^{n:k}(V) \ra \in 
[\M_g^{k+}(V) \wedge V^{n+}, B{\mathbb T}^{n+} \wedge \MU_{\Q}]_*$$ which 
record the tangent lines at the marked points; summing over $n$-fold cap 
products with a cycle ${\bf z} \in H_*(V,H^*(B{\mathbb T}))$ defines a theory 
based on configuration spaces with marked points restricted to lie on $\bf z$.
This recovers the Gromov-Witten potential and the WDVV family of quantum 
multiplications. \medskip

\noindent 
{\bf 2.5ii} Here is a sketch of a symplectic analogue, restricted for 
simplicity to surfaces with 
at most one boundary component. Let $\AC$
be the (contractible) space of almost-complex structures compatible with
the symplectic structure on $(V,\omega)$. A 
smooth map $u$ from a closed $\Sigma$ to $V$ will be called {\bf symplectic} 
if $u^*\omega$ is nondegenerate (and hence symplectic) on $\Sigma$; such 
a map induces a 
monomorphism of tangent bundles. Let $\proto(\Sigma;V)$ be the manifold 
of these maps, and let $\proto^*(\Sigma;V)$ be the bundle over
$\proto(\Sigma,V) \times \AC$ with fiber at $(u,J_{V})$
the (contractible) space of almost-complex structures on $\Sigma$
compatible with $u^*\omega. \; \Diff(\Sigma)$ acts freely on
this space, with quotient the bundle of spaces of
pseudoholomorphic maps from $\Sigma$ to $V$ over $\AC$. More generally, 
if $D \subset \Sigma$ is a smooth disk in the surface, and $x$ is
a point in its interior, let $\Diff(\Sigma,D,x)$ be the group of
orientation-preserving diffeomorphisms of $\Sigma$ which 
preserve both $D$ and $x$. There is then a bundle $\proto^*(\Sigma,
D,x;V)$ with fiber the space of almost-complex structures on $\Sigma - 
\partial D$ compatible with $u^*\omega$: to be precise, we want 
almost-complex structures on both sides, which extend continuously
to $\partial D$, but which need not agree there. Restriction to $D$ 
is a Fredholm map from $\proto^*(\Sigma,D,x;V)$ to $$\AC 
\times_{\Diff(\Sigma,D,x)} \widetilde {\el V}\;,$$ the last term in the 
product being the universal cover of the free loopspace of $V$. The range of 
this map is homotopy-equivalent to the product of $B\Diff(\Sigma)$ with the 
Borel construction [with respect to the action of the circle by rotating 
loops] on $\widetilde {\el V}$, and its index [11 Ch.\ 8] is $(1-g) 
{\Dim}_{\C}(V) + u^*(c_{1}(V))[\Sigma]$. The pullback of its cobordism class  
to the formal normal bundle to the fixed-point set is a kind of 
Gromov-Witten invariant in $MU^*(B\Diff(\Sigma)^+ \wedge {\bf {\hat H}F}(V))$, 
where ${\bf {\hat H}F}(V)$ is a pro-spectrum [as in the appendix to [3]] 
constructed from the normal bundle of the fixed loops. This is a conjectural 
analogue, for the Floer homotopy type associated to the area functional
on a free loopspace, of the theory sketched in \S 2.3. \medskip

\noindent
{\bf 2.6} Finally, I want to mention that the space $H^1(\Sigma,G)$ of
principle bundles [with compact Lie structure group $G$] defines a potential
Gromov-Witten map $$H^1(\Sigma,G) \times_{\Diff(\Sigma)}E\Diff(\Sigma)
\rightarrow H^1(\partial \Sigma,G) \times B\Diff(\Sigma)$$ which [as 
far as I know] has not yet been made the basis for a theory of topological
gravity \dots

\section{General nonsense}

\noindent
Monoidal functors between monoidal categories form a monoid, much as
homomorphisms between abelian groups form an abelian group. Manin 
and Zograf [7] suggest that we think of these families of theories as 
parametrized by the Picard group of invertible objects. Such objects 
can be identified with the points of ${\rm Spec} \; MU_{\Q}$, which 
are one-dimensional formal group laws; but there is no natural way to 
compose them. Kontsevich-Witten theory [4] suggests the natural parametrizing
object is the Hopf algebra $\mathcal Q$ of Schur $Q$-functions: these 
are Hall-Littlewood symmetric functions of the eigenvalues of a 
positive-definite matrix $\Lambda$, evaluated at $t = -1$. \medskip

\noindent
The Kontsevich-Witten genus $MU \rightarrow \mathcal Q$ defines 
a formal group law with $Q$-function coefficients; its exponential (aside
from normalization) is the asymptotic expansion as $\Lambda \rightarrow 
+\infty$ of the Mittag-Leffler exponential $$\sum_{n \geq 0} \frac{{\rm Tr}
\;{\Lambda}^n}{\Gamma (1 + \half n)} \; \;.$$ There is a natural 
{\bf twisted} charge one action of the Virasoro algebra on the $Q$-functions,
and the image under $\tau_{\bf kw}$ of the fundamental class 
$$\exp(\sum [\M_g]) \in H_*({\rm SP}^{\infty}(\coprod \M_g),\Q)$$ of the 
space of not-necessarily-connected curves is an $sl_2$-invariant 
highest-weight vector, or {\bf vacuum} state. \medskip

\noindent
Recently Eguchi {\it et al}, Dubrovin, Getzler [5], and others have begun
to extend this fundamental result of Kontsevich-Witten theory to topological
gravity coupled to quantum cohomology. The relevant Virasoro
representation appears to be defined on a group of loops on the torus 
$H^*(V,\R/\Z)$, twisted by the endomorphism $\half H + tX$, 
where $H,X,Y$ generate the standard $sl_2$ action on the Hodge cohomology
of $V$, and $t = c_1(V)/\omega$ depends on the K\"ahler class. \medskip

\noindent
This twisted torus is mysterious even when $V$ is a point. The adjoint
operation on $(\rm Gravity)_{1+1}$ defines an involution on the Picard group
of invertible theories, and the analogy with Abel-Jacobi theory suggests
that $\mathcal Q$ represents its skew-adjoint part. [The Hopf algebra 
defined by the cohomology of $BU$ represents the Witt ring functor,
defined on commutative rings by $$A \mapsto {\W}(A) = (1 + tA[[t]])^{\times}
\;;$$ away from two, the involution sending $a(t) \in \W(A)$ to $a(-t)$ 
splits $\W$ into eigenspaces. The cohomology of $\rm SO/U$ is the negative 
eigenspace.]

\newpage

\bibliographystyle{amsplain}

\end{document}